\documentclass[12pt]{article}
\usepackage[a4paper, total={6in, 8in}]{geometry}
\usepackage{xcolor}
\usepackage{algorithm}
\usepackage{algpseudocode}
\usepackage{graphicx}
\usepackage{hyperref}

\def\NN{{\cal N}}

\newtheorem{conjecture}{Conjecture}

\begin{document}

\title{\bf A fresh look to a randomized \\ massively parallel graph coloring algorithm} 

\author{Boštjan Gabrovšek \& Janez Žerovnik}

\newcommand{\Addresses}{{
\bigskip
\footnotesize

B.~Gabrovšek, \textsc{Univerza v Ljubljani, Fakulteta za strojništvo,
Aškerčeva cesta 6, 1000 Ljubljana, Slovenija;
Institute of Mathematics, Physics and Mechanics,
Jadranska ulica 19, 1000 Ljubljana, Slovenia.
}
\textit{E-mail address}: \texttt{bostjan.gabrovsek@fs.uni-lj.si}

\medskip

J.~Žerovnik , \textsc{Univerza v Ljubljani, Fakulteta za strojništvo,
Aškerčeva cesta 6, 1000 Ljubljana, Slovenija;
Institute of Mathematics, Physics and Mechanics,
Jadranska ulica 19, 1000 Ljubljana, Slovenia.
}
\textit{E-mail address}, M.~Dane: \texttt{janez.zerovnik@fs.uni-lj.si}

}}

\date{\today}

\maketitle

\begin{abstract}

Petford and Welsh introduced a sequential heuristic algorithm for (approximately) solving the NP-hard graph coloring problem. The algorithm is based on the antivoter model and mimics the behaviour of a physical process based on a multi-particle system of statistical mechanics. It was later shown that the algorithm can be implemented in massively parallel model of computation. The increase of processing power in recent years allow us to perform an extensive analysis of the algorithms on a larger scale,  leading to possibility of a more comprehensive understanding of the behaviour of the algorithm including the phase transition phenomena.
\end{abstract}


\section{Introduction}

Graph coloring is among the most studied NP-hard problems in combinatorial optimization. 
It is generally believed that NP-hard problems cannot be
solved to optimality within times which are polynomially bounded
functions of the input size. Therefore there is much interest in
heuristic algorithms which can find near-optimal solutions within
reasonable running times for a large proportion of instances. 

Based on the antivoter model of Donnely and Welsh \cite{DonellyWelsh},
Petford and Welsh proposed a randomized algorithm for graph coloring,
which performs very well on some random graphs \cite{PeWe,DM}, as well
as on some other types of graphs \cite{Jann}. 
With some straightforward modifications, the algorithm was also very competitive on 
frequency assignment problems \cite{rovinj, speedup,cejore}.

Petford and Welsh experimented with their algorithm on a model of random 3-colorable graphs.
They observed that there are some combinations of parameters of random graphs that are extremely hard for their algorithm, 
while otherwise their algorithm was running on average in linear time. 
Having no theoretical eplanation, Petrod and Welsh write that the curious behavior ``is not unlike the phenomenon of phase-transition 
which occurs in the Ising model, Potts model and other models of statistical mechanics''.

We have recently designed an algorithm for clustering that is motivated by this coloring algorithm. 
The results were very promising \cite{Ikica}. 
This has motivated us to better understand the basic algorithm, which we do by experimental study that is reported here.
We test the algorithm on a larger amount of data in order to have a better understanding of its performance (number of parallel steps).
In particular, we run and analyse the following experiments:
\begin{itemize}
 \item in (\cite{parcomp,corr}) the algorithms performs poorly on a couple of instances. We run several experiments and test the hypothesis that the algoritem performance depends only on the average number of degree for random or regular graphs,
 \item in order to test whether the performance of the algorithm depends on the average degree, we test performance by varying the average degree of the graph and the temperature (parameter $T$) of the system. We test the performance on   random graphs and on $r$-regular graphs (both known to be colorable),
 \item we test whether the algorithm behaves similarly for higher degree colorings ($k=4,5,\dots, 10$) and confirm the existence of critical regions that are related to the phase transition phenomena.
\end{itemize}

The rest of paper is organized as follows. 
In Section \ref{sec:alg} we recall the $k$-coloring decision problem and the 
algorithms described in \cite{parcomp} and \cite{corr}. 
In section \ref{sec:exp} we present the new experiments and analyse the results.

\section{The problem and the algorithm} \label{sec:alg}

\bigskip
{\bf Graph coloring problem.}
The $k$-coloring decision problem is a well known NP-complete problem
for $k\geq 3$. It reads as follows:

\medskip
\fbox{\em
\begin{minipage}{8.5cm}
Input: graph $G$, integer $k$ \\
Question: is there a proper $k$-coloring of $G$?
\end{minipage}}
\medskip


We say a mapping $c: V(G) \to \NN$ is a {\em proper coloring} of $G$ if it assigns
different colors to adjacent vertices. 
Any mapping $c: V(G) \to \NN$ will be called a {\em coloring}, 
and will be considered as one of the feasible solutions 
of the problem. 
We introduce the cost function
$E(c)$ 
to be the number of {\em bad} edges, i.\ e.\ edges with both ends colored by the same color by coloring $c$. 
Proper colorings are exactly the colorings for which $E(c)=0$ and finding a coloring $c$ with $E(c) =0$ is 
equivalent to answering the above decision problem where the coloirng constructed is a {\em witness} $c$ proving the correctness of the answer.

\bigskip
{\bf The algorithm of Petford and Welsh.}
The basic algorithm \cite{PeWe} starts with a random initial 3 coloring of the input graph and then applies an iterative process. 
In each iteration a vertex creating a conflict is chosen at random. 
The random distribution is uniform among the bad vertices and the chosen vertex is recolored according to some probability distribution. 
The color distribution favours colors which are less represented in the neighborhood of the chosen vertex, see the 
expression (\ref{BasicFormulae}) below.
The algorithm has a straightforward generalization to $k$ coloring 
(taking $k=3$ gives the original algorithm) \cite{DM}. 

In a pseudo language the algorithm of Petford and Welsh can be written as


\begin{algorithm}
 \caption{Petford-Welsh algorithm}\label{alg:pw}
\begin{algorithmic}[1]
 \State color vertices randomly with colors $1,2,\dots,k$
 \While{not stopping condition} 
 \State select a bad vertex $v$ (randomly)
 \State assign a new color to $v$
 \EndWhile
\end{algorithmic}
\end{algorithm}

Bad vertex is selected uniformly random among vertices which are endpoints of some bad (e.g. monocromatic) edge.
A new color is assigned at random. The new color is taken from the set
$ \{1, 2, \ldots, k\} $. 
Sampling is done according to probability distribution defined as follows:

The probability $p_i$ of color $i$ to be chosen as a new color of vertex $v$ 
is proportional to
\begin{equation}
 p_i \approx \exp( -S_i /T), \label{BasicFormulae}
\end{equation}
where $S_i$ is the number of edges with one endpoint at $v$ and
with color $i$ assigned to the other endpoint. 
Petford and Welsh used $ 4^{-S_i}$ which ie equivalent to using $T \simeq 0. 72$ in (\ref{BasicFormulae}).
(Because $ \exp( -x/T) = 4^{-x}$ implies $T \simeq 0. 72$.)

The stopping criteria is either reaching a time limit (in terms of the number of calls to the function which computes a new color)
or if a proper coloring is found. 
If no proper coloring is found, the solution with minimal cost $E(c)$
is reported and can be regarded as an approximate solution to the problem. 
However, there is no guarantee on the quality of the solution known. 

\bigskip
{\bf Connection to simulated annealing and the generalized Boltzmann machine.}
Here we briefly discuss the parameter $T$ of the algorithm. 
 $T$ is the parameter of the algorithm, which may be called temperature because of the
analogy to the temperature of the simulated annealing algorithm and to the temperature of the generalized Boltzmann machine neural network as observed in \cite{Ann}. 
With term {\em simulated annealing} we refer to the optimization heuristics as proposed, for example, in \cite{KiGV83}.
Generalized Botzmann machine is a generalization of the
Boltzmann machine, a popular neural network that is based on a stochastic spin-glass model 
and widely used in artificial intelligence \cite{Hinton}.
The main difference is that the generalized Boltzmann machine as defined in \cite{Ann} uses
multistate neurons in contrast to bipolar neurons of the usual Boltzmann machine.

These analogies are based on the following simple observation. 
Pick a vertex and denote the old color of the chosen vertex by $j$ and the new color by $i$. 
The number of bad edges $E^\prime$ after the move is
$$ E^\prime = E - S_j + S_i
$$
where $E$ is the number of bad edges before the change. 
We define 
$\Delta E = E -E^\prime = S_j -S_i$. 
At each step, $j$ is fixed and hence $S_j$ and $E$ are fixed. 
Consequently, it is equivalent to define 
the probability of choosing the new color $i$ to be proportional to either
 $\exp ( -S_i /T)$,
$ \exp( \Delta E /T)$ or
 $\exp ( E^\prime /T)$. 

In order to see the relation to the Boltzmann machine,
 recall that the number of bad edges is a usual definition of energy function both in simulated annealing and in 
generalized Boltzmann machine with multistate neurons. 
Therefore, the algorithm PW is in close relationship 
to constant temperature operation of the generalized Boltzmann machine (for details, see \cite{Ann} and the references there). 
The major difference is in the firing rule. 
While in Boltzman machine all neurons are fired with equal probability, in algorithm 
of Petford and Welsh, only bad vertices are activated. 

As already explained, the original algorithm of Petford and Welsh 
uses probabilities proportional to $4^{-S_i}$, which corresponds to $T \approx 0. 72$. 
Larger values of $T$ result in higher probability of accepting a move which increases the number of bad edges. 
Clearly, a very high $T$ results in chaotic behavior similar to a pure random walk among the colorings ignoring the 
their energy. 
On the other hand, with low values of $T$, the algorithm behaves very much like iterative improvement, 
quickly converging to a local minimum.

There is another slight difference with the usual implementation of simulated annealing heuristics (SA) \cite{KiGV83} as an optimization algorithm. Namely, in SA, a change improving the cost is always accepted while in the other case, the acceptance probability 
is used which depends both on the difference of costs and the temperature. 
In algorithm of Petford and Welsh, all changes are taken according to the probability using (\ref{BasicFormulae}).
The cost improving changes thus may not be accepted, although this happens very rarely in the majority of cases.

Thus, our algorithm is similar to both simulated annealing heuristics and to sequential operation of the Boltzmann machine.
Its acceptance probabilty is (at given $T$) practically equivalent to Boltzmann machine.
As the Boltzmann machine is a highly parallel asynchronous device, a comparison with parallel implementations 
is even more interesting. 
Here, we recall the parallel version of the algorithm of Petford and Welsh that 
differs from the generalized Boltzmann machine only in the firing rules in both phases of its operation.

\bigskip
{\bf Parallel algorithm.}
In \cite{parcomp},  a massively parallel version of Algoritm \ref{alg:pw} was proposed.
The naive algorithm (Algorithm \ref{alg:pwpNAiVE}) was later improved in \cite{corr} by a version that runs in two phases 
thus avoiding the looping that may appear at some configurations within the instance.
The improved algorithm (Algorithm \ref{alg:pwpar}) 
first aims to find a $2k$ coloring and does not recolor all bad vertices simultaneously because each 
change is only done with some probability (i.\ e.\ 0.6). 
In the second phase, the result of the first phase provides independent sets of vertices which may be recolored in parallel 
without any conflict.
The resulting algorithm stil shows maximal speedup in comparison to the original version, e.\ g.\
the instances solved in linear time sequentially are expected to be solved in constant time in parallel \cite{corr}.
The algorithm formaly reads as Algorithm \ref{alg:pwpar}.

Note that in the first phase, the firing rule is nearly equivalent to the same of Boltzmann machine in which each neuron (vertex) wakes up at some (random) time and performs recoloring. Note that there is no synchronization among the neurons.
While sinhronization is not of particular importance in the first phase of our algorithm, in the second phase, it is essential that synchronization based on the result of the first phase is used.

\begin{algorithm}[htb]
 \caption{Massively parallel variant of the Petford Welsh algorithm (naive version)}\label{alg:pwpNAiVE}
 \begin{algorithmic}[1]
 \State color vertices randomly with colors $1,2,\ldots,k$
 \While{not stopping condition}
 \State $bad\_vertices \gets \{ v \mid v \mbox{ is bad} \}$ 
 \For{\textbf{all} $v \in bad\_vertices$}
 \State assign a new color to $v$
 \EndFor
 \EndWhile 
 \hspace{13em}\raisebox{2.1\baselineskip}[0pt][0pt]{$\left.\rule{0pt}{1.6\baselineskip}\right\}\ \mbox{in parallel}$}
 \end{algorithmic}
\end{algorithm}

\begin{algorithm}[hhtb]
 \caption{Massively parallel variant of the Petford Welsh algorithm}\label{alg:pwpar} 
 \begin{algorithmic}[1]
 \Procedure{MPPW\_phase1}{G}
 \State color vertices randomly with colors $1,2,\ldots,k, k+1,\ldots, 2k$
 \While{not stopping condition}
 \State $bad\_vertices \gets \{ v \mid v \mbox{ is bad} \}$ 
 \For{\textbf{all} $v \in bad\_vertices$}
 \State assign a new color to $v$ with a probability of $60\%$
 \EndFor
 \EndWhile 
 \\
 \Return{coloring}
 \EndProcedure
 \hspace{24em}\raisebox{3.7\baselineskip}[0pt][0pt]{$\left.\rule{0pt}{1.5\baselineskip}\right\}\ \mbox{in parallel}$}
 
 \bigskip
 \Procedure{MPPW\_phase2}{$G$}
 \State $c \gets $ \Call{MPPW\_phase1}{G}
 \State color vertices randomly with colors $1,2,\ldots, k$
 \While{not stopping condition}
 \State $bad\_vertices \gets \{ v \mid v \mbox{ is bad} \}$ 
 \For{\textbf{all} $v \in bad\_vertices$}
 \If{$step \bmod (2k) = c(v)$}
 \State assign a new color to $v$
 \EndIf
 \EndFor
 \EndWhile 
 \\
 \Return{coloring}
 \EndProcedure
 \hspace{14em}\raisebox{6.0\baselineskip}[0pt][0pt]{$\left.\rule{0pt}{1.5\baselineskip}\right\}\ \mbox{in parallel}$}
 \end{algorithmic} 
\end{algorithm}
 
\section{Experiments}\label{sec:exp}

In our experiments, we use two classes of graphs. The first class are graphs of the form $$G(n,k,p),$$
where $n$ is the number of vertices, $k$ is the number of partitions and $p$ is the probability two
vertices from distinct partitions are adjacent (see \cite{PeWe,parcomp,code}).

The second class are $d$-regular graphs of the form $$R(n,k,d),$$
where $n$ is the number of vertices, $k$ is the number of partitions and $d$ is the degree of vertices.
The python code for generating such graphs can be found in \cite{code}. In short, 
the algorithm splits vertices into $k$ partitions and connects, in random order, each point to $d$ random points in distinct partitions. 
If there are points left, that are not of degree $d$ and cannot be connected, we delete an edge and add two other edges, in particular, if $u$ and $v$ are vertices with $\deg(u) < d$ and $\deg(v) < d$, then we find an edge $u'v'$, such that $u$ and $u'$ do not belong to the same partition and $v$ and $v$ do not belong to the same partition. We delete the edge $u'v'$ and add edges $uu'$ and $vv'$ to the graph.

In both cases, partitions are of equal size if
$k$ divides $n$, otherwise their sizes differ by at most one vertex.

\begin{figure}[ht]\centering
 \includegraphics[scale=1.0]{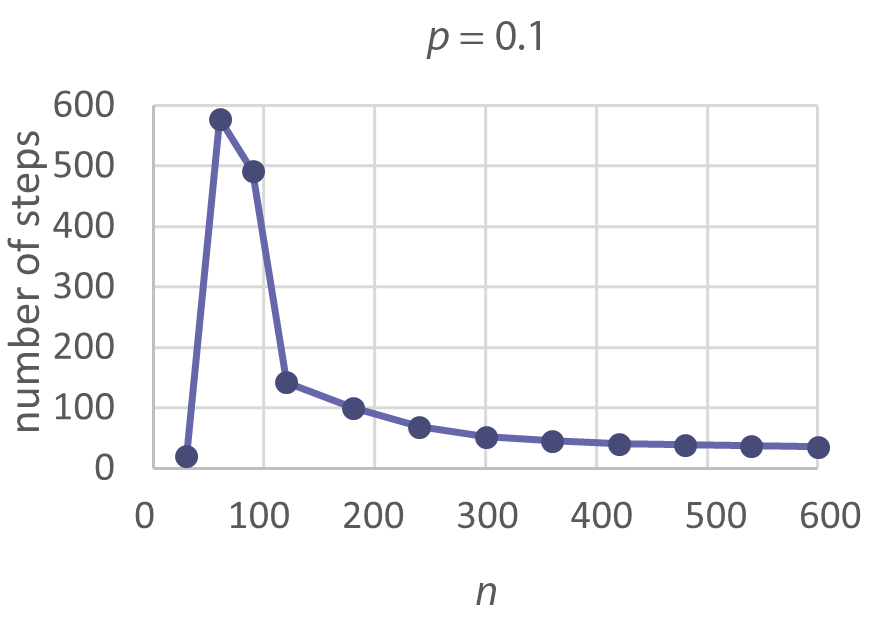} \qquad \includegraphics[scale=1.0]{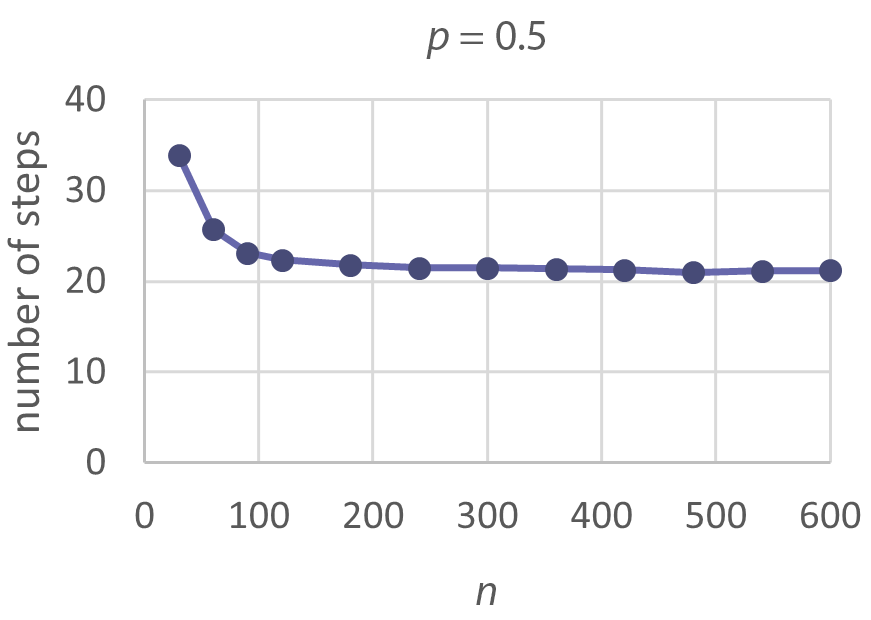}
 \caption{Performance of IMPPW (parallel steps as a function of $n$), cf.~Table~1 in \cite{parcomp}.}\label{diagram:old}
\end{figure}

\bigskip
{\bf Preliminary experiment.}
With code in \cite{code} we reproduced the results from~\cite{parcomp, corr}, with much larger sample size, $n=10000$ (instead of $n=100$). The results are presented in Figure~\ref{diagram:old}, 
confirming that our implementation runs exactly the same algorithm as the original.

Observe that the algorithm does not perform well only on a narrow interval, which we call the
{\em critical region}. 
Observing some experimental results, 
limited by the computing resources available at the time, the following conjecture was proposed \cite{DM}.

\begin{conjecture}\label{conj}
The critical regions are characterized by the equation
\begin{equation}
    \frac{2pn}{k} \approx \frac{16}{3}.\label{critical}  
\end{equation}
\end{conjecture}

This conjecture generalizes the conjecture of Petford and Welsh who observed that the equation 
 $\frac{2pn}{3} \approx \frac{16}{3}$ is valid within the critical region \cite{PeWe}.
More precisely, they observed that given $p$, the graphs $G(n,3,p)$ with $n \approx \frac{8}{p}$ are likely hard instances for the algorithm.

After we confirm the basic observations in the main references, 
we continue with experimental results that may shed some more light to the behaviour of the algorithm and possibly to some more 
general phenomena. 
In particular, we wanted to understand better, if and how the hard instances can be related to the average degree of the graphs. 
In relation to this, we wish to check whether the conjecture above captures the main information that determines the critical regions.
Furthermore, we are interested in question 
``what is the effect of the temperature'' (or, equivalently, the basis of exponent expression (\ref{BasicFormulae}))
on the performance of the algorithm.
Below we provide results of the experiments with some comments that answer some and, at the same time, open some new questions.

\bigskip
{\bf First experiment.}
In the first experiment, we show that, with a fixed number of components, the critical region depends on the average degree $\overline{d}(G)$ of the graph, where 
$$\overline{d}(G) = \frac{n p (k-1)}{k}.$$

\begin{figure}[hhtb]\centering
 \includegraphics[scale=0.7]{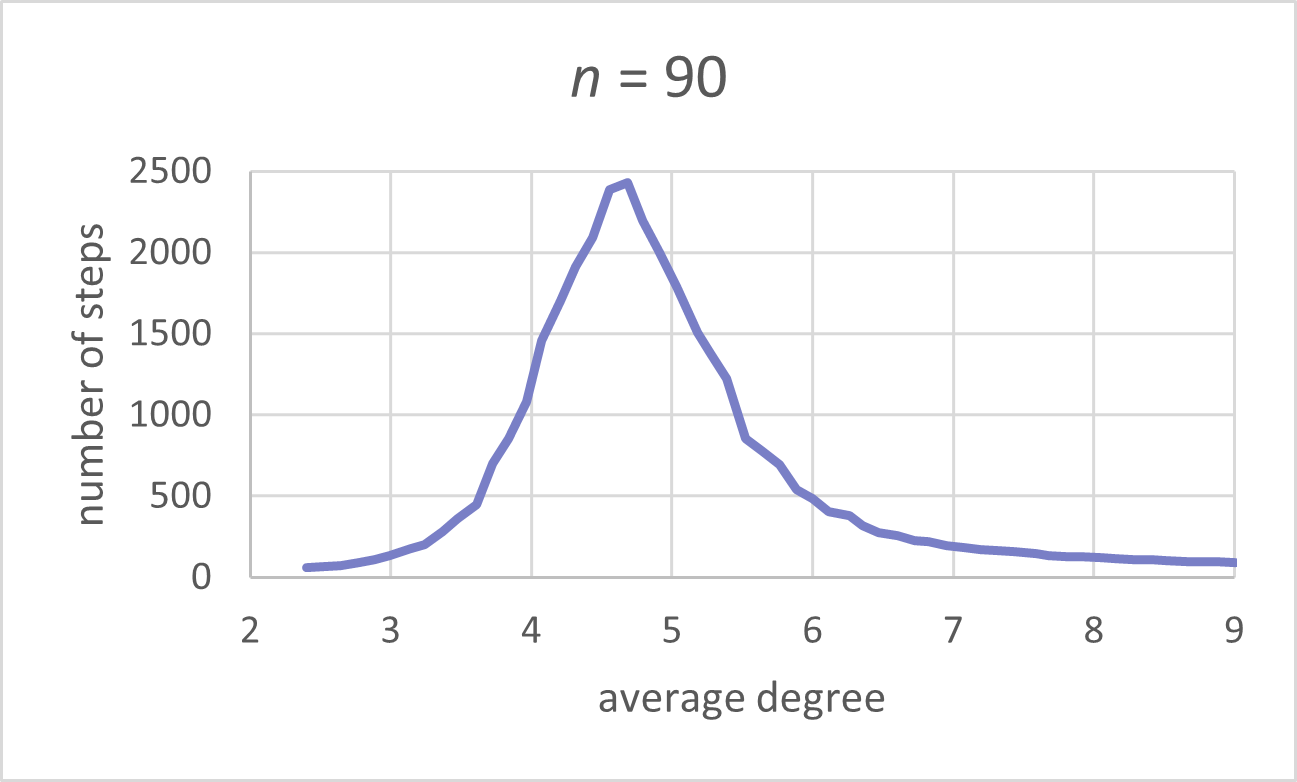} \qquad \includegraphics[scale=0.7]{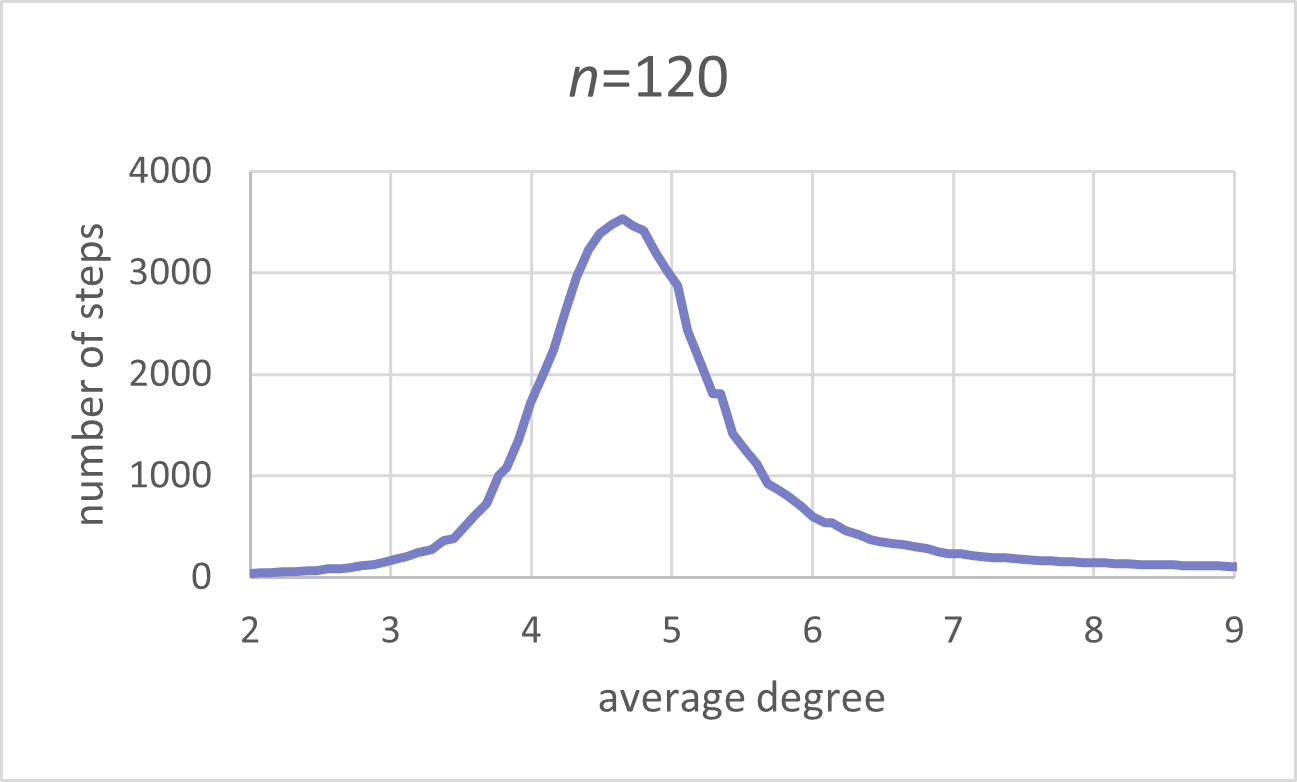}\\
 \includegraphics[scale=0.7]{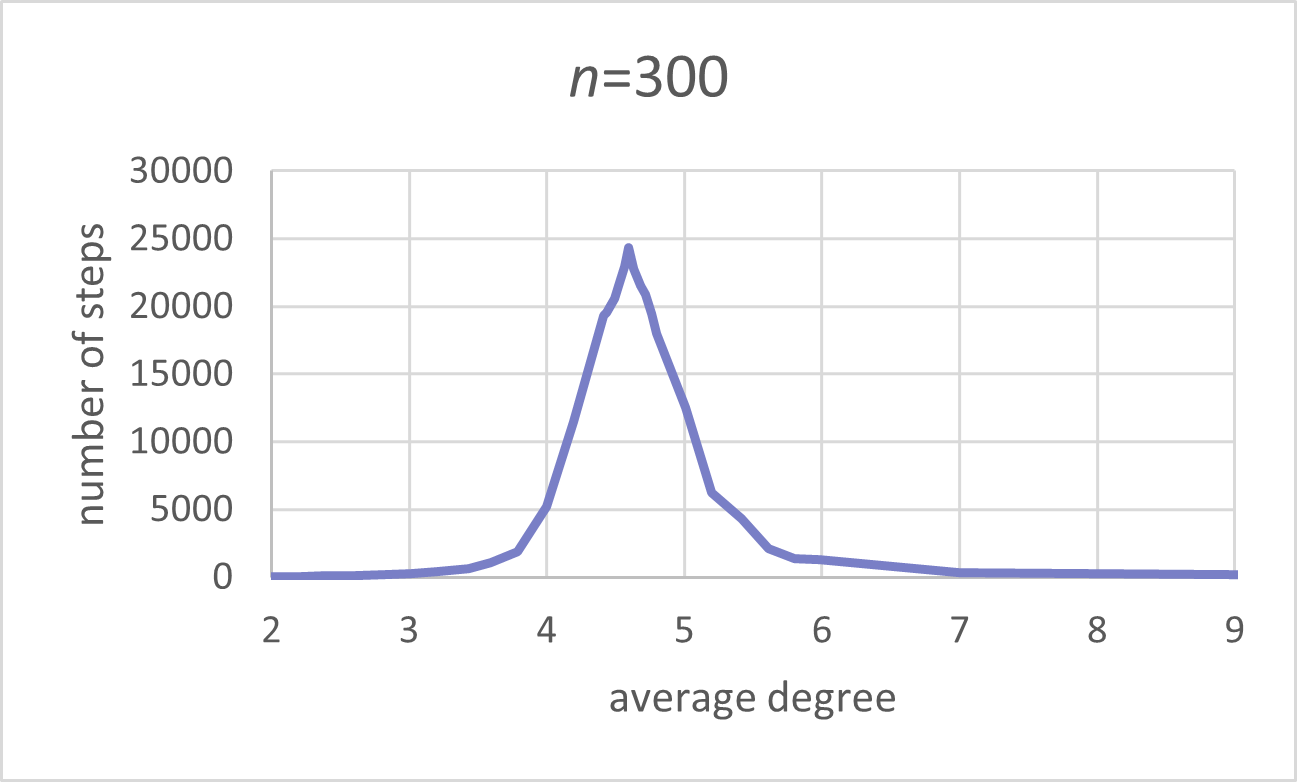} \qquad \includegraphics[scale=0.7]{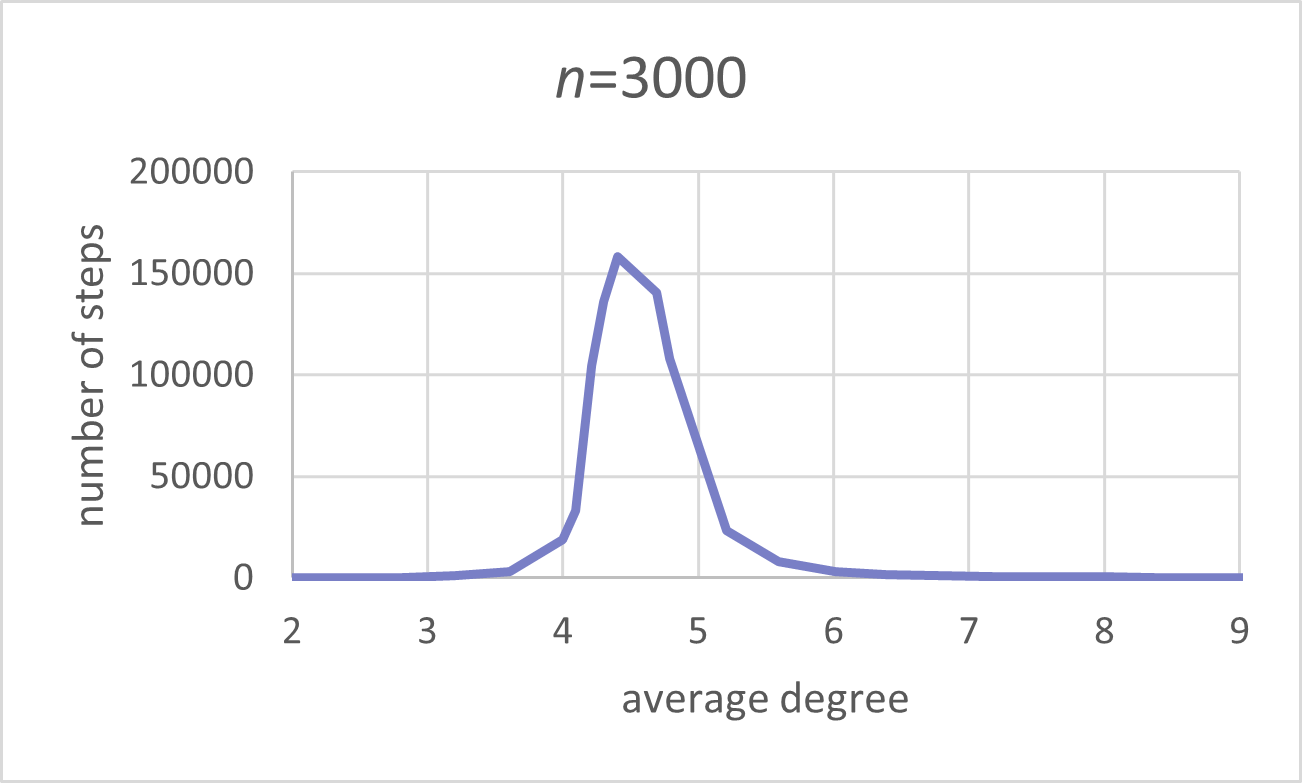}
 \caption{Performance of IMPPW (steps as a function of average degree)}\label{diagram:ex2}
\end{figure}

We choose four datasets: graphs of classes $G(90,3,p)$, $G(120,3,p)$, $G(300,3,p)$, and $G(3000,3,p)$. 
The sample size is $10000$ for $n \in \{90, 120, 300\}$ and 
a bit smaller, $2000$, for $n = 3000.$
We choose the parameter $p$ in such a way that the average degrees of each class varies from 2 to 9.
The results are presented in diagrams on Figure~\ref{diagram:ex2}.

Note that in all four experiments, the extreme value does not depend on the number of vertices.
It is obvious that the hard instances are in the interval where the average degree is between 4 and 5, even more, the peaks are 
 within the range 4.5 -- 4.7 in all four diagrams.

The phenomena can be explained as follows.
Observe that before the critical region ($\overline d(G) < 4.5$), the partitions (e.g. sets of vertices of the same color in a proper coloring)
 are loosely defined and there are multiple optimal solutions (see supplementary video 1).
On the other hand, after the critical region ($\overline d(G) > 4.7$), the algorithm converges fast, since the partitions are densly connected and thus very well defined (see supplementary video 3). For a graph in the critical region see supplementary video 2. All three videos can be also accessed at~\cite{code}. 

\bigskip
 
\begin{figure}[hhbbt]\centering
 \includegraphics[scale=0.75]{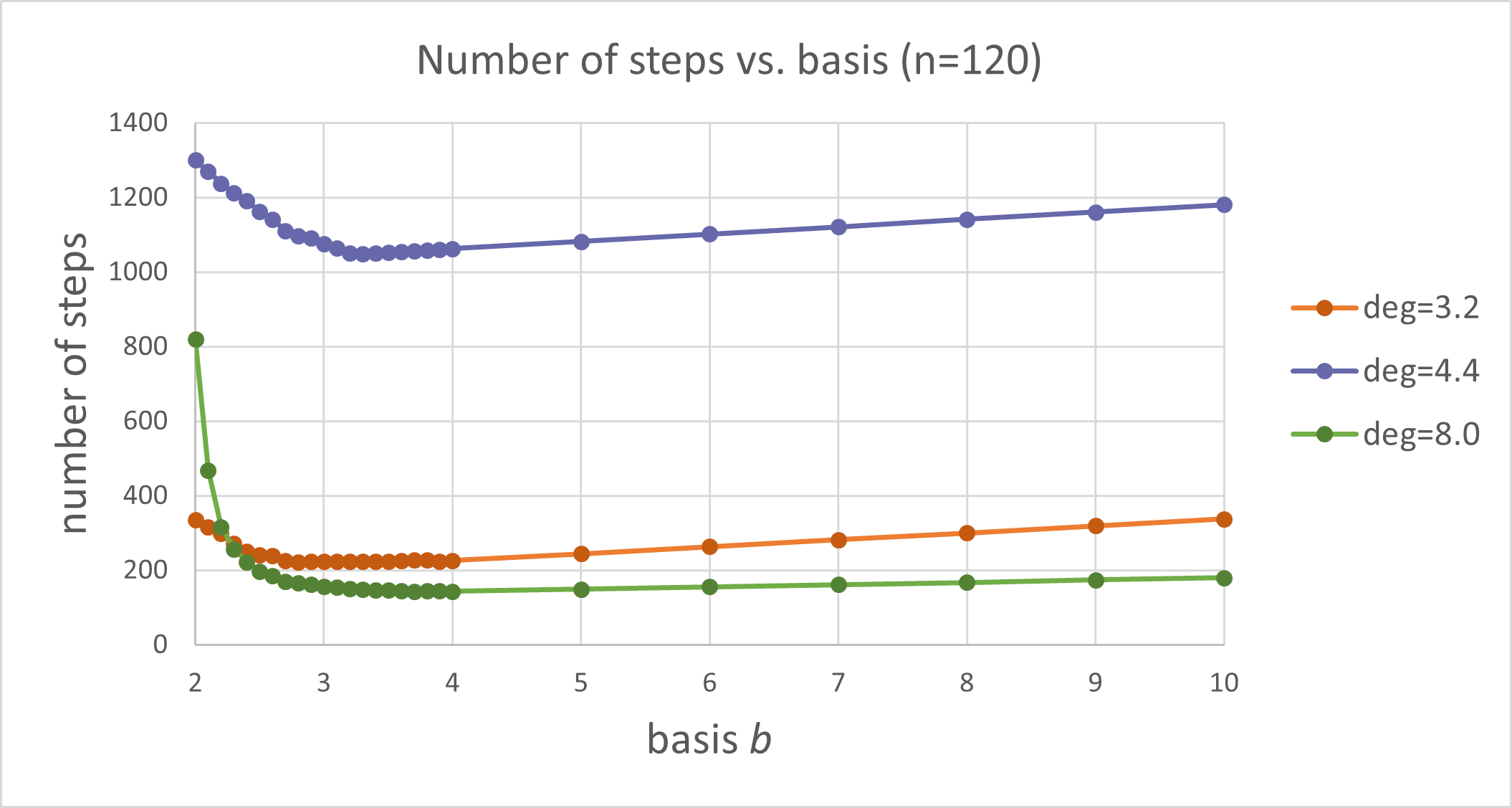}
 \caption{Performance of IMPPW (basis vs. number of steps)}\label{diagram:ex3a}
\end{figure}

\begin{figure}[hhbbt]\centering
 \includegraphics[width=0.7\linewidth]{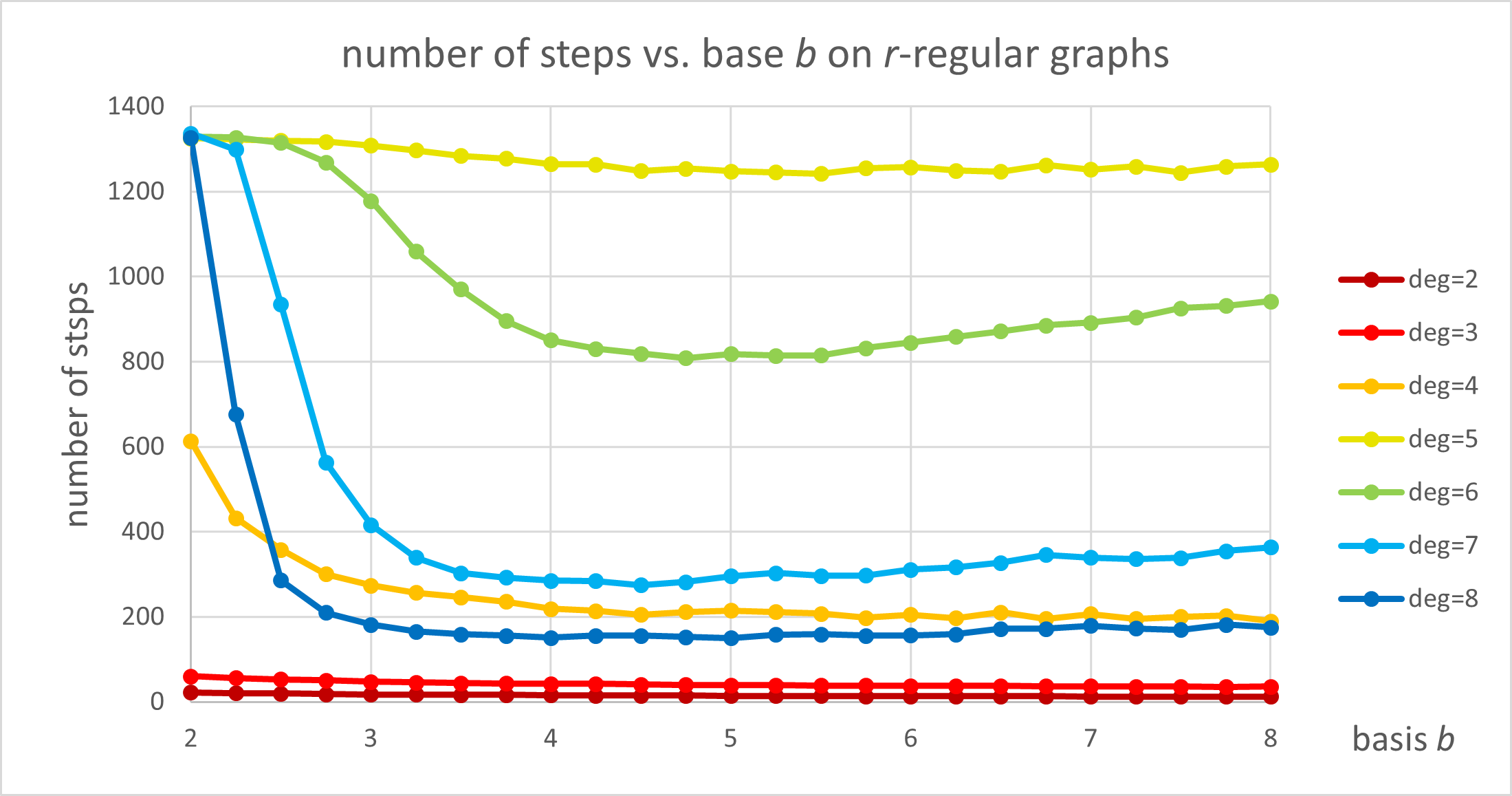}
 \caption{Performance of IMPPW (basis vs. number of steps for $d$-regular graphs)}\label{diagram:ex3b}
\end{figure}

{\bf Second experiment.}
In this experiment, we vary the parameter $b$ over the values from $2$ to $10$ and measure the performance of the algorithm in the classes $G(120,3,p)$, where $p$ is chosen in such a way that the average degrees are 3.2, 4.4, and 8.0. 
The results are presented in diagram~\ref{diagram:ex3a}.

We conclude that the basis $b$ does not dramatically influence the performance. 
Seemingly, the values between $3.0$ and $4.0$ are well behaved and indeed, taking any value $b \in [3, 4]$ performes similarly.
We repeat the experiment on the graphs $R(120,3,d)$, where $d=2, 3, \ldots, 8$.
Results are presented in diagram on Figure~\ref{diagram:ex3b} and clearly confirm the earlier observations.

The question ``what is the optimal value $b$?'' may be expected to have an answer between 3 and 4. 
Note that this is a question equivalent to the question which is the optimal temperature of the simulated annealing 
(i.\ e.\ ``annealing'' with constant temperature) which, however, seems to be rather complex problem \cite{Cohn,icannga}.
(Recall that   $ \exp( -x/T) = b^{-x}$ so $b= \exp(1/T)$.)
As our insight is limited by the special class of instances used, we do not wish to dig further in the question of optimal $b$. 
On the other hand, we can confirm that, likely, the algorithm is rather robust regarding the choice of $b$ (and equivalently) to the choice of parameter $T$).

\begin{figure}[hhtb]\centering
    \includegraphics[width=0.48\linewidth]{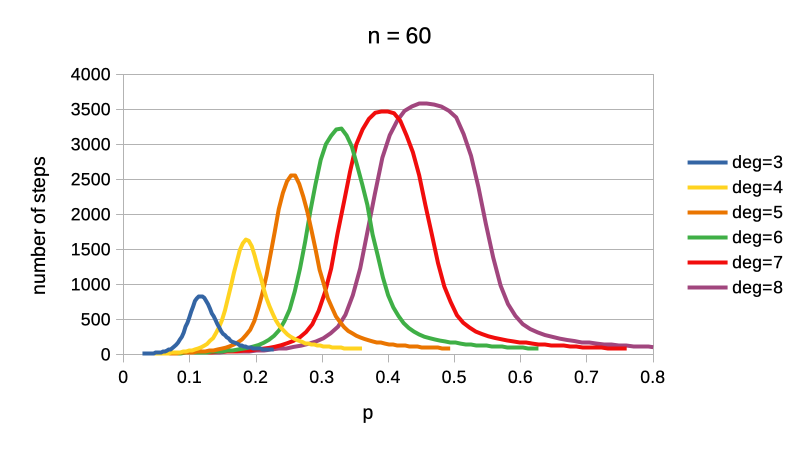}
    \includegraphics[width=0.48\linewidth]{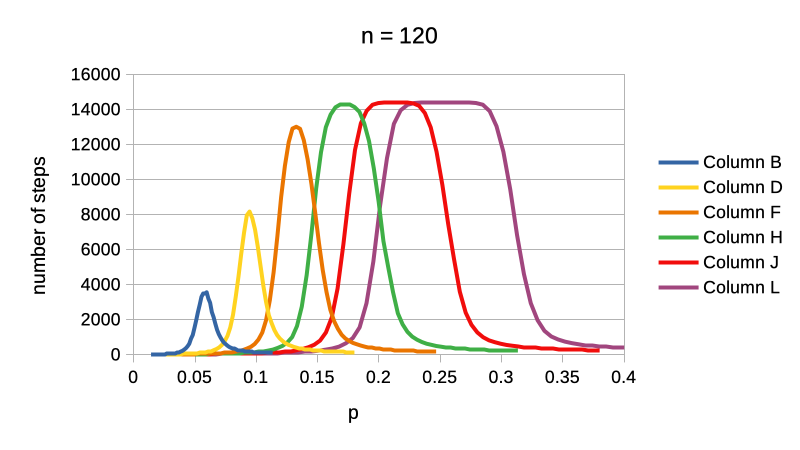}
    \includegraphics[width=0.48\linewidth]{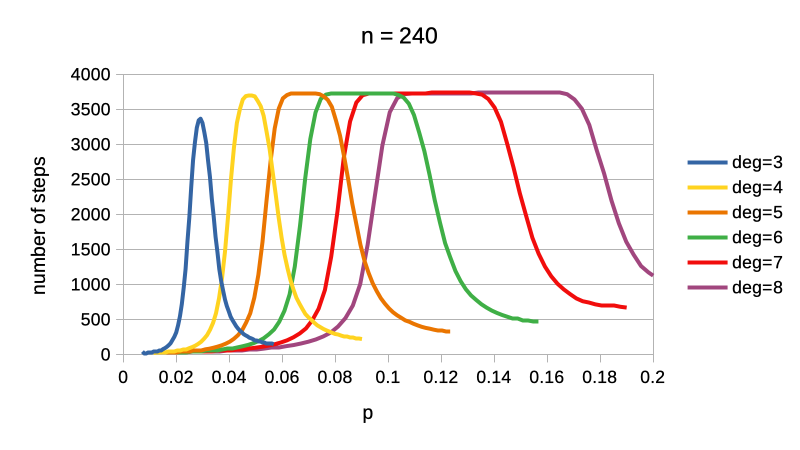}
    \caption{Performance of IMPPW for $k$-colorings. The graphs show how the number of steps depend on the parameter $p$.
    }
    \label{diagram:a}
   \end{figure}

\bigskip
{\bf Third experiment.}
In this experiment, we compute run times for graphs in $G(n,k,p)$
where we fix $k = 3, 4, \ldots, 8$ partitions and $n = 60, 120, 240$ vertices
and in each case observe how the number of steps needed depends on the probability $p$.
 In all the experiments, the sample size is $5000$.
The results are presented in Figure \ref{diagram:a}.

\begin{figure}[hhtb]\centering
    \includegraphics[width=0.48\linewidth]{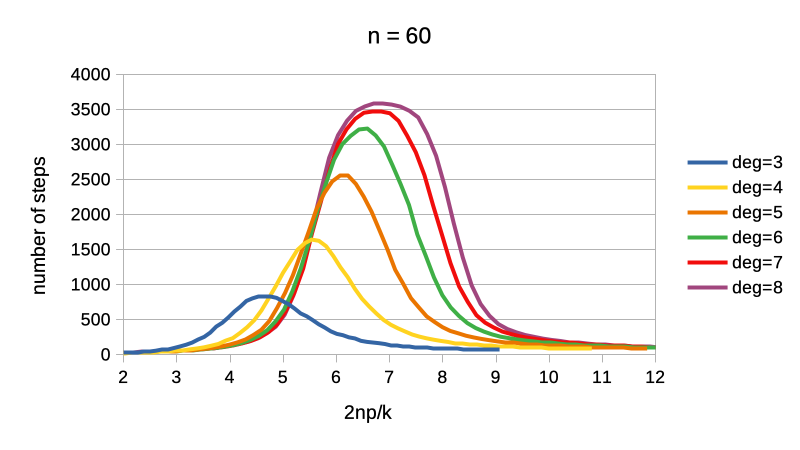}
    \includegraphics[width=0.48\linewidth]{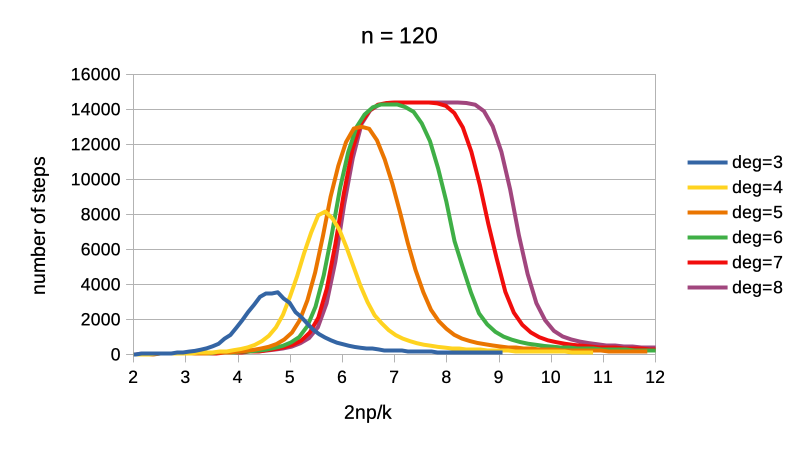}
    \includegraphics[width=0.48\linewidth]{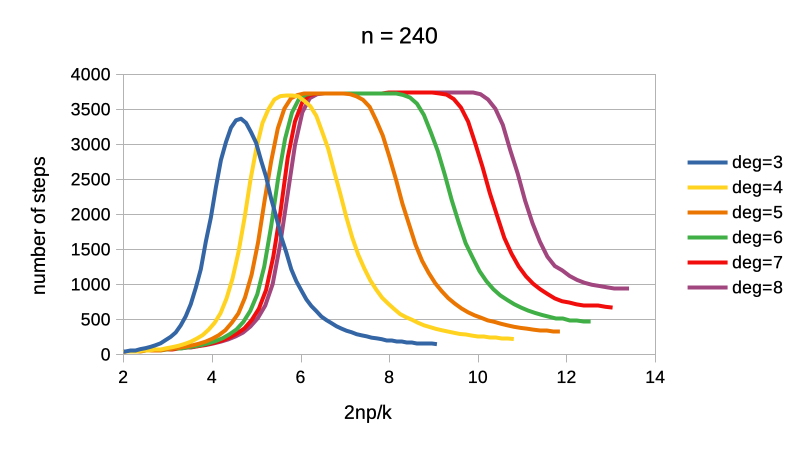}
    \caption{Performance of IMPPW for $k$-colorings. The graphs show how the number of steps depend on the parameter $p$.
    }
    \label{diagram:b}
   \end{figure}

According to Conjecture \ref{conj}, the hard instances are characterized by some constant value of $\frac{2pn}{k}$.
However, in Figure \ref{diagram:b} we observe how the performance depends on $\frac{2pn}{k}$.
We conclude that the critical region is not characterized exactly by $\frac{2pn}{k}$, thus Conjecture \ref{conj} should be replaced by a  better one. 

\begin{figure}[hhtb]\centering
    \includegraphics[width=0.48\linewidth]{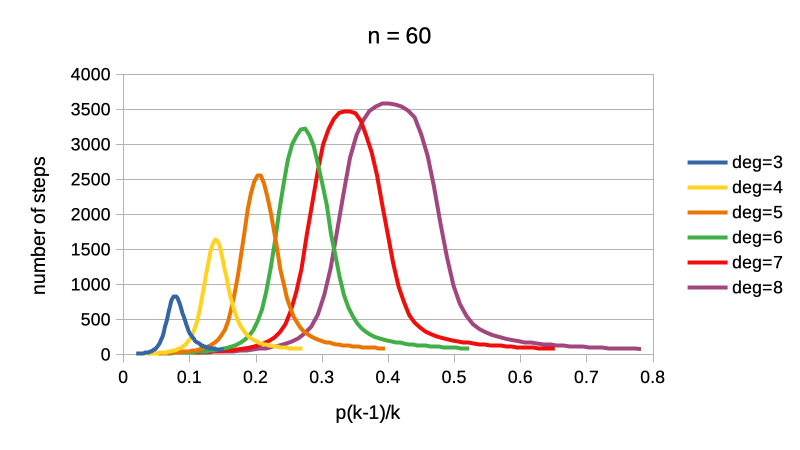}
    \includegraphics[width=0.48\linewidth]{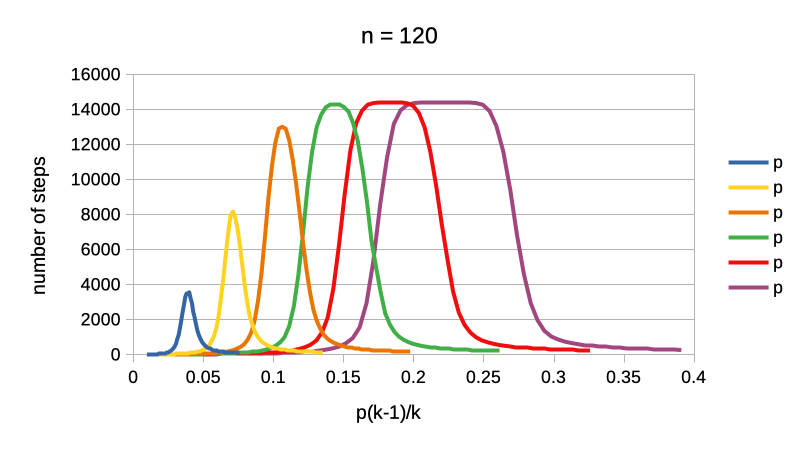}
    \includegraphics[width=0.48\linewidth]{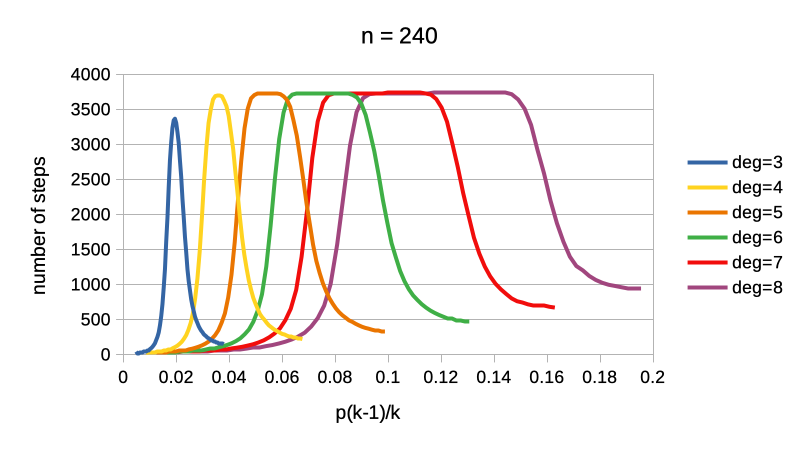}
    \caption{Performance of IMPPW for $k$-colorings. The graphs show how the number of steps depend on the average degree, which is proportional to $\frac{p(k-1)}{k}$.
    }
    \label{diagram:d}
   \end{figure}

In Figure~\ref{diagram:d} we plot how the number of steps depend on the expression $\frac{p(k-1)}{k}$, which is 
proportional, if $n$ is fixed, to the average degree. 
We observe from the figure that the critical region can not be explained only by the average degree.

\begin{figure}[hhtb]\centering
    \includegraphics[width=0.48\linewidth]{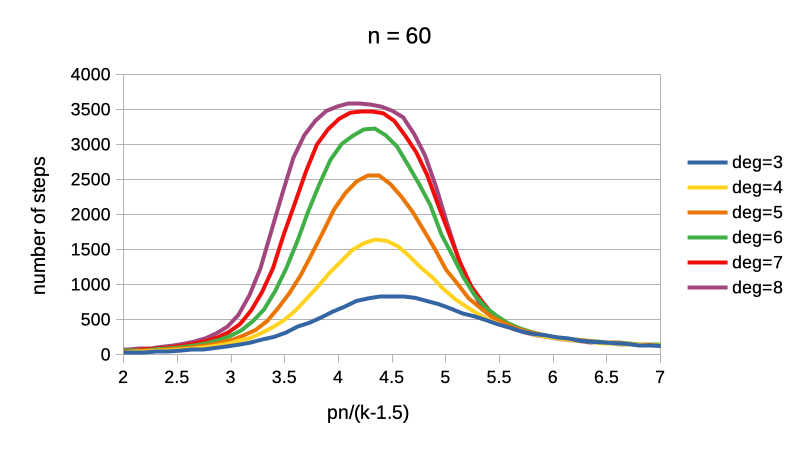}
    \includegraphics[width=0.48\linewidth]{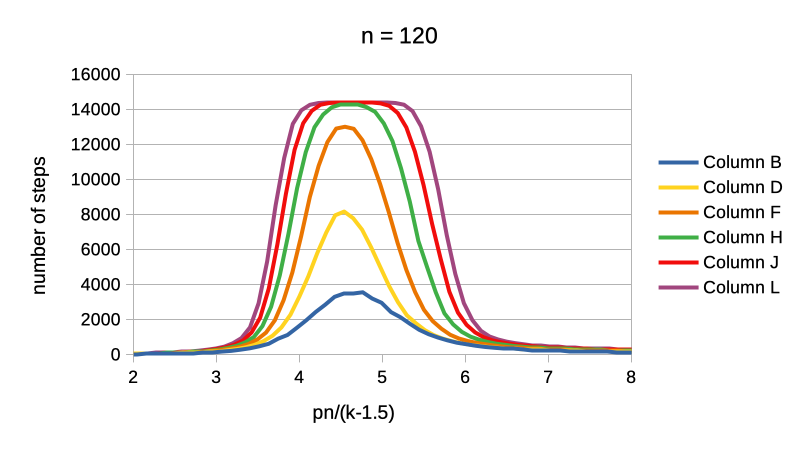}
    \includegraphics[width=0.48\linewidth]{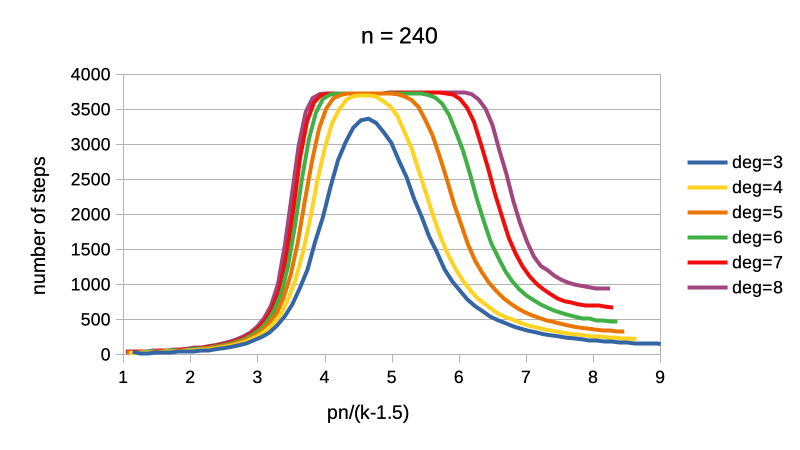}
    \caption{Performance of IMPPW for $k$-colorings. The graphs show how the number of steps depend on the expression $\frac{pn}{k-1.5}$. The cut-off is due to the time constraint of the algorithm.}
    \label{diagram:c}
   \end{figure}

At present, we do not have a good idea how to improve the conjecture to better characterize the critical region with an 
expression that would have some natural meaning. 
We have tested some slight modifications, and found that the expression $\frac{p}{k-1.5}$ remains fairly constant for varied $k$.
See Figure~\ref{diagram:c}, where we plot how the number of steps depend on the expression $\frac{p}{k-1.5}$.

Therefore, it seems that the critical region is approximately characterized by the equation
$$\frac{np}{k-1.5} \approx 4.3.$$

\section{Conclusions}\label{sec:conclusions}

Parallel version of Petford and Welsh $k$-coloring algorithm was extensively tested on two classes of random graphs. 
The conjecture on existence of a critical region where the algorithm has nearly prohibitively long run times 
is confirmed for the case $k=3$, while a generalized conjecture \cite{DM} is shown to need an adjustment.
More precisely, 
we have observed that the critical region appears where $\frac{p}{k-1.3}$ holds. 
This result opens at least two interesting questions. 

\begin{itemize}
\item 
can the property $\frac{np}{k-1.5} \approx 4.3$ be naturally explained as some feature of the instances ?
\item 
is there another expresion that fits the data, and has some meaning which can explain the structure of hard instances ?
\end{itemize}

Treshold phenomena have attracted a lot of attention in the context of random combinatorial problems \cite{editorial}
and in theoretical physics \cite{Philathong2021, j}. 
In statistical physics, phase transitions have been studied for more than a century. Let us only mention here that the spin glasses, 
a purely theoretical concept, has triggered a new branch of theoretical physics that resulted in a Nobel prize attributed to Giorgio Parisi in 2022
\cite{Mezard}.
Graph coloring with $k=3$ colors has been considered in several papers, see \cite{CulbGent2001,Boettcher2004} and the references there. 
In contrast to graph classes of $k$ colorable graphs used in this paper, the usual random graph model considered are 
graphs $G(n,p)$ where each of the possible $\frac{n(n-1)}{2}$ edges appears independently with probability $p$.
Not surprisingly, for 3-coloring, it is found that the critical mean degree where the phase transition occurs is 
around $\alpha_{crit} \approx 4.7$, for example the estimate 4.703 was put forward  in \cite{Boettcher2004}.
In the same paper, the analysis implies that the hardest instances are among graphs with average degree between 4.42 and $\alpha_{crit}$.
It seems that the $k$-coloring was not considered, hence we can not compare our findings about the hardest instances with previous work.
We conclude that 
further study of the critical regions is a promising avenue of research that may have some implications that go beyond only deeper 
understanding of behavior of the algorithm of Petford and Welsh.

\section*{Acknowledgements}\label{sec:ack}
The authors were supported in part by the Slovenian Research Agency,
 grants J2-2512, J1-4031, N1-0278 (Boštjan Gabrovšek), and J2-2512 and P2-0248 (Janez Žerovnik).


\def\by{ }
\def\jour{ }
\def\vol{ \bf }
\def\pages{ }
\def\paper{ \it }
\def\book{ \it }
\def\yr{ }

\Addresses

\end{document}